\documentclass{article}
\usepackage{amsfonts,amsthm,latexsym,amsmath,amssymb,amstext,amsxtra,amscd,epsfig,psfrag}
\usepackage{graphics,graphicx, bezier, float, color}
 \usepackage{url}
\usepackage[hmargin=1.2in,vmargin=1.5in]{geometry}
 \usepackage[table]{xcolor}
\usepackage{caption}
\usepackage{relsize}
\usepackage[pdfpagelabels,draft,implicit=false]{hyperref}
 \usepackage{calc}
\usepackage{bm}
\usepackage[utf8]{inputenc}
\newtheorem{theorem}         {Theorem}
\newtheorem{prop+}           {Proposition}
\newtheorem{coro+}           {Corollary}
\newtheorem{lemm+}           {Lemma}

\theoremstyle{definition}

\newtheorem{not+}            {Notation}
\newtheorem{example}         {Example}[section]

\newtheorem{question}         {Question}
\theoremstyle{remark}
\newtheorem{rema+}{\bf Remark}

 \newenvironment{proposition}{\begin{prop+}}{\end{prop+}}
 \newenvironment{corollary}{\begin{coro+}}{\end{coro+}}
 
 \newenvironment{remark}{\begin{rema+}}{\end{rema+}}

\theoremstyle{definition}
\newtheorem{definition}{Definition}

\newenvironment{prf}[0]{\textbf{Proof:}}{\hfill$\Box$}

\newcommand{\bq}{\begin{eqnarray}}
\newcommand{\eq}{\end{eqnarray}}
\newcommand{\beq}{\begin{eqnarray*}}
\newcommand{\eeq}{\end{eqnarray*}}

\date{}
\begin{document}\numberwithin{equation}{section}

\title{\textbf{$\text{Red}$-injective modules}}
 \author{Juma Kasozi, David Ssevviiri\footnote{Corresponding author}  and Vincent Umutabazi}
\maketitle
\begin{center}
 Department of Mathematics\\
 Makerere University, P.O BOX 7062, Kampala Uganda\\
 E-mail addresses:  kasozi@cns.mak.ac.ug, ssevviiri@cns.mak.ac.ug, umutabazivincent@yahoo.fr
\end{center}
\begin{abstract}
  Let $\text{Red}(M)$ be the sum of all reduced submodules of a module $M$.
  For modules over   commutative rings, $\text{Soc}(M)\subseteq \text{Red}(M)$.
  By drawing motivation from how $\text{Soc}$-injective modules were defined by Amin et. al.  in 
   \cite{amin2005},  we introduce $\text{Red}$-injective modules,
   study their properties and use them to characterize   quasi-Frobenius rings and $V$-rings.
  \end{abstract}
  \textbf{Keywords}: Injective modules; $\text{Red}$-injective modules;
  $\text{Soc}$-injective modules; Quasi-Frobenius rings; $V$-rings

  \textbf{MSC 2010} Mathematics Subject Classification: 16D50, 16D60, 16L60
\maketitle

\section{Introduction}

  For a not necessarily commutative ring $R$, Lee and Zhou  in \cite{lee2004reduced} defined an
  $R$-module $M$ to be reduced if for all $r\in R$ and  $m\in M$, $mr=0$
  implies that $Mr\cap mR=\{0\}$.   This definition is equivalent to saying that for all $r\in R$ and
  $m\in M$, $mr^{2}=0$ implies that $mRr=\{0\}$, see \cite{ssevviiri2013contribution} for the proof.  However, for modules over commutative rings we get
  Definition \ref{def} below.

\begin{definition}\label{def}
  An $R$-module $M$ is {\it reduced} if for all $r\in R$ and $m\in M$,
  $mr^{2}=0$ implies that $mr=0$.
\end{definition}

 Except in Example \ref{ex}, all rings are unital, commutative and associative.
  Modules are right unital defined over rings.
  A submodule is reduced if it is reduced as a module. A submodule of
  a reduced module is reduced but a factor module of a reduced module
  need not be reduced.
  The $\mathbb{Z}$-module $\mathbb{Z}$ is reduced but its factor module
  $\mathbb{Z}/n\mathbb{Z}$ is not reduced for a non-square free integer
  $n$.
  The socle of an $R$-module $M$, denoted by $\text{Soc}(M)$ is the sum
  of simple submodules of $M$.
  Let $\text{Red}(M)$ denote the sum of reduced submodules of $M$, i.e.,
  $$\text{Red}(M):=\sum_{i\in I}\{N_{i}~|~N_{i}~\text{is a reduced submodule of}~ M\}.$$

\begin{definition}
  An $R$-module $M$ is {\it semi-reduced} if $\text{Red}(M)=M$.
\end{definition}

\begin{proposition}\label{prop1}
  For any $R$-module $M$, the following implications hold:
  $$\text{simple}\Rightarrow \text{semi-simple}\Rightarrow \text{reduced}\Rightarrow \text{semi-reduced}.$$ 
\end{proposition}

\begin{prf}
  We prove that a semi-simple module is reduced. The other implications
  follow from the definition of semi-simple and semi-reduced modules
  respectively.
  Since a simple module is prime\footnote{An $R$-module $M$ for which $RM\neq \{0\}$
 is {\it prime} if for all $a\in R$ and every $m\in M$, $am=0$ implies that $m=0$ or $aM=\{0\}$.}
 and every prime module is reduced, a simple module is reduced. Suppose that $M$ is a
 semi-simple module and $mr^{2}=0$ where $m\in M$ and $r\in R$.
 Then, $(m_{1}, m_{2}, \cdot \cdot \cdot, m_{i}, \cdot \cdot \cdot)r^{2}=0$
 where $(m_{1}, m_{2}, \cdot \cdot \cdot , m_{i}, \cdot \cdot \cdot)=
 m \in M=\bigoplus_{i\in I}M_{i}$ for some simple modules $M_{i}$.
 Since every simple module is reduced,  $m_{i}r^{2}=
 0\Rightarrow m_{i}r=0~\forall ~i\in I$. Hence,  $mr=0$, and $M$ is reduced.
 \end{prf}

 \begin{corollary}
  For any $R$-module $M$, $\text{Soc}(M)\subseteq \text{Red}(M)$. 
\end{corollary}

\begin{prf}
  The proof follows from the fact that a semi-simple module is semi-reduced
  which is proved in Proposition \ref{prop1}.
\end{prf}

  Note that for semi-simple modules and for modules without nonzero reduced
  submodules, $\text{Soc}(M)=\text{Red}(M)$.
\begin{example}
  A reduced module need not be semi-simple.
  $\mathbb{Z}$ and $\mathbb{Q}$ are reduced $\mathbb{Z}$-modules but they
  are not semi-simple.
\end{example}


\subsection{Other basic definitions}
\begin{definition}{\rm \cite[Definition 1.1]{amin2005}}
 Let $M$ and $N$ be $R$-modules. $M$ is \emph{socle-$N$-injective}
 ($\text{Soc}$-$N$-injective) if any $R$-homomorphism
 $f: \text{Soc}(N)\rightarrow M$ extends to $N$. Equivalently, for any
 semi-simple submodule $K$ of $N$, any $R$-homomorphism
 $f : K \rightarrow M$ extends to $N$. An $R$-module $M$ is
 {\it $\text{Soc}$-quasi-injective} if $M$ is $\text{Soc}$-$M$-injective.
 $M$ is {\it $\text{Soc}$-injective} if $M$ is $\text{Soc}$-$R$-injective.
 $R$ is right (self-) {\it $\text{Soc}$-injective}, if the module
 $R_{R}$ is $\text{Soc}$-injective (equivalently, if $R_{R}$ is
 $\text{Soc}$-quasi-injective).
\end{definition}

\begin{definition}{\rm \cite[Definition 1.2]{amin2005}}
 An $R$-module $M$ is called \emph{strongly $\text{Soc}$-injective}, if $M$ is
 $\text{Soc}$-$N$-injective for all $R$-modules $N$. A ring $R$ is called
 {\it strongly $\text{Soc}$-injective}, if the module $R_{R}$ is
 strongly $\text{Soc}$-injective.
\end{definition}
  Definitions 3 and 4 together with Corollary 1 motivate us to have Definitions 5
  and 6 respectively.
\begin{definition}
 An $R$-module $M$ is called \emph{$\text{Red}$-$N$-injective} if any $R$-homomorphism
 $f:K\rightarrow M$ extends to $N$ for any semi-reduced submodule $K$ of $N$.
 $M$ is called {\it $\text{Red}$-quasi-injective} if it is $\text{Red}$-$M$-injective.
 $M$ is called {\it $\text{Red}$-injective} if it is $\text{Red}$-$R$-injective.
\end{definition}
\begin{definition}
 An $R$-module $M$ is called \emph{strongly-$\text{Red}$-injective}, if $M$ is
 $\text{Red}$-$N$-injective for all $R$-modules $N$.
\end{definition}

 In Definition \ref{de}, we recall  different generalizations of injective modules that we
 later use in the sequel. As with $\text{Soc}$-injective and $\text{Red}$-injective
 modules defined above, these generalizations of injective modules were  defined by
 relaxing conditions on the lifting property of homomorphisms.

 \begin{definition}\label{de}
 If $M$ and $N$ are $R$-modules, then
 \begin{enumerate}
   \item $M$ is {\it $N$-injective} if every $R$-homomorphism from a submodule of $N$ into $M$ can be
   extended to an $R$-homomorphism from $N$ into $M$.
   
   \item $M$ is {\it quasi-injective} if it is $M$-injective.
    \item  $M$ is {\it $N$-simple-injective} if for any submodule $L$ of $N$,    any homomorphism $\theta: L \rightarrow M$ with $\theta (L)$ simple,
    can be extended to a homomorphism $\beta : N \rightarrow M$.
   \item $M$ is {\it simple-injective} if it is simple $R$-injective.
   \item $M$ is {\it strongly simple-injective}, if $M$ is simple-$N$-injective
   for all right $R$-modules  $N$.
   \item $M$ is {\it min-$N$-injective} if, for every simple submodule $L$ of $N$,
   every homomorphism $\gamma : L \rightarrow M$ extends to $N$.
   \item $M$ is {\it min-injective} if it is min-$R$-injective.
   \item $M$ is {\it strongly min-injective}, if it is min-$N$-injective for all
   $R$-modules $N$.
   
   \item $M$  is {\it pseudo-injective}  if any monomorphism from a submodule of $M$ to $M$ extends to an endomorphism of $M$.
 \end{enumerate}
\end{definition}

\subsection{Notation}
Throughout this paper, $N\subseteq^{e} M$, $N\oplus M$, $N\subseteq^{\oplus} M$,
 and $N\leq M$, mean that $N$ is an essential submodule of $M$, a direct sum of $N$
 and $M$, $N$ is a direct summand of $M$, $N$ is a submodule of $M$ respectively.

\subsection{Paper  roadmap}

In Section 1, we have given the introduction, defined key terms,  given the notation used  and the roadmap for the paper.

 Section 2 is devoted to obtaining  properties of Red-injective modules and their generalizations. An equivalent definition of a Red-injective
 module is obtained. It is shown that any injective module is strongly Red-injective and a Red-injective module is Soc-injective. Other 
 implications with known generalizations of injective modules are given. 
The class of (strongly) $\text{Red}$-injective  $R$-modules is closed under isomorphisms, direct products 
 and summands. If $M$ is a Noetherian module, then a direct sum of Red-$M$-injective is Red-$M$-injective. For a family of $R$-modules
 $\{M_i~:~i\in I\}$, an $R$-module $N$  is Red-$(\oplus_{i\in I}M_i)$-injective if and only if it is Red-$M_i$-injective for each $i$.
  For a  projective $R$-module $M$,  every quotient of a $\text{Red}$-$M$-injective 
 $R$-module is $\text{Red}$-$M$-injective if and only if  $\text{Red}(M)$ is projective if and only if every quotient of an injective $R$-module
  is Red-$M$-injective.   Over a principal ideal domain a free module is $\text{Red}$-injective if each of its
 submodule is $\text{Red}$-injective.  Red$(N)$-lifting modules are introduced. It is shown  that if a module $N$  is Red($N$)-lifting, then 
 any $R$-module $K$ is Red-$N$-lifting if and only if $K$ is $N$-injective. It is shown that  
 $\text{Red}$-quasi-injective modules inherit  a weaker version of C2-condition and C3-condition.

 In Section 3, we characterize quasi-Frobenius rings and right $V$-rings
 in terms of strongly $\text{Red}$-injective modules.   A ring $R$ is quasi-Frobenius if and only if every strongly Red-injective $R$-module is projective.
 A ring $R$ is a right $V$-ring if and only if every simple $R$-module is strongly Red-injective. A question is raised as to whether  
  $\text{Red}$-quasi-injective modules  and  $\text{Soc}$-quasi-injective
 modules  are clean and or satisfy the exchange property.

\section{$\text{Red}$-injective modules}
\begin{proposition}
 For $R$-modules $K, ~M$ and $N$, the following statements are equivalent:
 \begin{enumerate}
   \item Any $R$-homomorphism $f:K\rightarrow M$ extends to $N$ for any semi-reduced          submodule $K$ of $N$.
   \item Any $R$-homomorphism $f:\text{Red}(N)\rightarrow M$ extends to $N$.
 \end{enumerate}

\end{proposition}
\begin{prf}
 \begin{description}
   \item[$1\Rightarrow 2$] since $\text{Red}(N)$ is semi-reduced.
   \item[$2\Rightarrow 1.$] Suppose $f:K\rightarrow M$ is an $R$-homomorphism
        and $K$ is a semi-reduced submodule of $N$. Since $K\leq \text{Red}(N)$, then
        $f$ extends to $N$.
 \end{description}
\end{prf}

\begin{proposition}

  If $N$ is an $R$-module, then 
\begin{enumerate}
 \item any injective module is strongly $\text{Red}$-injective,
 \item a $\text{Red}(N)$-injective module is $\text{Soc}(N)$-injective.
\end{enumerate} 
\end{proposition}
 
\begin{prf}
 \begin{enumerate}
   \item Let $M$ be an injective module. Then $M$ is $N$-injective for
         every $R$-module $N$. For every submodule $K$ of $N$, any
         $R$-homomorphism $f: K\rightarrow M$ extends to $N$.
         For every module $N$, any $R$-homomorphism
         $f:\text{Red}(N)\rightarrow M$ extends to $N$. Hence, $M$
         is strongly $\text{Red}$-injective.
   \item Suppose $f:\text{Soc}(N)\rightarrow M$ is an $R$-homomorphism
         and $M$ is $\text{Red}(N)$-injective. By Proposition 1,
         $\text{Soc}(N)$ is a semi-reduced submodule of $N$. Hence by
         Definition 5, $f$ extends to $N$. Thus, $M$           is $\text{Soc}$-$N$-injective.
 \end{enumerate}
\end{prf}

  Every projective module over a  right Noetherian
  right self-injective ring is strongly $\text{Red}$-injective.
  Let $R$ be a ring for which each module $M$ has $\text{Red}(M)=\{0\}$.
  Then, $M$ is strongly $\text{Red}$-injective.

 A $\text{Red}$-injective module need not be injective.
 The module $\mathbb{Z}_{\mathbb{Z}}$ is $\text{Red}$-injective but not injective.

\begin{theorem}  Let $\{ M_{i}:i\in I\}$ be a family of $R$-modules and
  $N$, $M$, $A$, $C$, $S$ and $K$ be $R$-modules. Then the following conditions hold:
\begin{enumerate}
  \item A direct product $\prod_{i\in I}M_{i}$  is $\text{Red}$-$N$-injective
        if and only if each $M_{i}$ is $\text{Red}$-$N$-injective.
  \item For $S\leq N$, if $M$ is
        $\text{Red}$-$N$-injective, then $M$ is $\text{Red}$-$S$-injective.
  \item For $M\cong N$; $M$ is $\text{Red}$-$S$-injective if and only
        if $N$ is $\text{Red}$-$S$-injective.
  \item For $A\cong B$; $C$ is $\text{Red}$-$A$-injective if and only if
        it is $\text{Red}$-$B$-injective.
  \item For $N\subseteq^{\bigoplus} M$,
        if $M$ is $\text{Red}$-$K$-injective, then $N$ is $\text{Red}$-$K$-injective.
\end{enumerate} 
\end{theorem}

\begin{prf}
\begin{enumerate}
\item We prove only for $M=M_{i}\times M_{j}$ where $i,~j\in I$.
      The proof for the general case is analogous. Let $M_{i}$ and $M_{j}$
      be $\text{Red}$-$N$-injective $R$-modules, $h:\text{Red}(N)\rightarrow N$
      and $f:\text{Red}(N)\rightarrow M_{i}\times M_{j}$ be any $R$-homomorphisms.\\
      Define $$f_{M_{i}}:\text{Red}(N)\rightarrow M_{i}~ \text{such that}~ \pi_{_{M_{i}}}\circ f=
      f_{M_{i}}$$ and $$f_{M_{j}}:\text{Red}(N)\rightarrow M_{j}~ \text{such that}~ \pi_{_{M_{j}}}\circ f=
      f_{M_{j}},$$ where $\pi_{_{M_{i}}}:M_{i}\times M_{j}\rightarrow M_{i}$
      and $\pi_{_{M_{j}}}:M_{i}\times M_{j}\rightarrow M_{j}$ are $R$-homomorphisms.
      Since $M_{i}$ and $M_{j}$ are $\text{Red}$-$N$-injective there exists
      $f'_{M_{i}}:N\rightarrow M_{i}$ and $f'_{M_{j}}:N\rightarrow M_{j}$
      such that $$f_{M_{i}}=f'_{M_{i}}\circ h~\text{and}~f_{M_{j}}=f'_{M_{j}}\circ h.$$
      By the uniqueness part of the universal property of direct product there exists
      an $R$-homomorphism $f':N\rightarrow M_{i}\times M_{j}$ such that $f=f'\circ h$.
      It follows that $\pi_{_{M_{i}}}\circ (f'\circ h)=f_{M_{i}}$
      and $\pi_{_{M_{j}}}\circ (f'\circ h)=f_{M_{j}}$.
      By the uniqueness of the universal property we conclude that $f=f'\circ h$.
      Hence, $f:\text{Red}(N)\rightarrow M_{i}\times M_{j}$ extends to $N$.
      Thus $M_{i}\times M_{j}$ is $\text{Red}$-$N$-injective.
      Conversely, assume that $M_{i}\times M_{j}$ is $\text{Red}$-$N$-injective.
      Let $h:\text{Red}(N)\rightarrow N$ and $f_{M_{i}}:\text{Red}(N)\rightarrow M_{i}$
      be any $R$-homomorphisms. Choose $f_{M_{j}}:\text{Red}(N)\rightarrow M_{j}$
      to be the zero $R$-homomorphism. We obtain $f':N\rightarrow M_{i}\times M_{j}$
      such that $f=f'\circ h$. Finally we obtain $f_{M_{i}}=\pi_{M_{i}}\circ f=(\pi_{M_{i}}\circ f')\circ h$.
      Hence $\pi_{M_{i}}\circ f':N\rightarrow M_{i}$ is an extension of $f_{M_{i}}$.
      Thus, $M_{i}$ is $\text{Red}$-$N$-injective. Similarly, $M_{j}$ is $\text{Red}$-$N$-injective.
  \item Consider the diagram in Figure 1, where $M$ is $\text{Red}$-$N$-injective.
  \begin{figure}[H]
  \centering
  \input rimwe.pic
  \label{rimwe}
  \caption{}
\end{figure}

  Since $S\leq N$, $\text{Red}(S)\leq \text{Red}(N)$. Consider inclusion maps
  $$k:\text{Red}(S)~\rightarrow ~\text{Red}(N);
  ~g:\text{Red}(S)~\rightarrow ~S;~h:\text{Red}(N)\rightarrow N~ \text{and}~ \iota:S~\rightarrow ~N.$$
  $f'\circ l:S\rightarrow M$ is an extension for any $R$-homomorphism
  $q:\text{Red}(S)\rightarrow M$. Thus $M$ is $\text{Red}$-$S$-injective.
\item
     Let $N\cong M$ where $\theta: N\rightarrow M$ is an $R$-isomorphism between them.
     Let $f_{N}:\text{Red}(S)\rightarrow N$ be any $R$-homomorphism. Since $M$ is
     $\text{Red}$-$S$-injective, any $R$-homomorphism $f_{M}:\text{Red}(S)\rightarrow M$
     extends to $f'_{M}:S\rightarrow M$. So for any $R$-homomorphism
     $h:\text{Red}(S)\rightarrow S$, $f_{M}=f'_{M}\circ h$.
     Since $M$ and $N$ are isomorphic there exists an inverse homomorphism
     $\theta^{-1}:M\rightarrow N$ such that $\theta^{-1}\circ f'_{M}:S\rightarrow N$
     is an $R$-homomorphism. Define $f'_{N}=\theta^{-1}\circ f'_{M}:S\rightarrow N$.
     Then, $f'_{N}$ is an extension of $f_{N}$.
     Thus, $N$ is $\text{Red}$-$S$-injective. Similarly, if $N$ is $\text{Red}$-$S$-injective
     then $M$ is $\text{Red}$-$S$-injective.
\item Suppose that $A\cong B$ and $C$ is $\text{Red}$-$A$-injective. We show that $C$
     is $\text{Red}$-$B$-injective.
\begin{figure}[H]
  \centering
  \input kabili.pic
  \label{kabili}
  \caption{}
\end{figure}

    Consider the diagram in Figure 2, where $f'_{A}:A\rightarrow C$
    is the extension of $f_{A}:\text{Red}(A)\rightarrow C$.
    Let also $f_{B}:\text{Red}(B)\rightarrow C$ be an $R$-homomorphism.
    Define $f'_{B}=f'_{A}\circ \theta :B\rightarrow C$.
    Then $f'_{B}:B\rightarrow C$ is the extension of $f_{B}$.
    Thus $C$ is $\text{Red}$-$B$-injective. A similar argument works
    for the converse.

  \item Let $N\subseteq^{\oplus} M$ and $M$ be $\text{Red}$-$K$-injective.
        We show that $N$ is $\text{Red}$-$K$-injective. Since $N\subseteq^{\oplus} M$,
        there exists an $R$-submodule $N'$ of $M$ such that $N\bigoplus N'=M$.
        Let $\pi_{_{N}}:N\bigoplus N'\rightarrow N$ be the projection $R$-homomorphism.
        Since $M$ is $\text{Red}$-$K$-injective, any $R$-homomorphism
        $f_{M}:\text{Red}(K)\rightarrow M$ extends to $f'_{M}:K\rightarrow M$.
        Suppose $f_{N}=\pi_{_{N}} \circ f_{M}:\text{Red}(K)\rightarrow N$.
        Define $f'_{N}=\pi_{_{N}} \circ f'_{M}:K\rightarrow N$.
        Then $f'_{N}:K\rightarrow N$
        is the extension of $f_{N}$. Hence, $N$ is $\text{Red}$-$K$-injective.
\end{enumerate}
\end{prf}
\begin{corollary}
          Let $N$ be an $R$-module, then 
   \begin{enumerate}
  \item a finite direct sum of $\text{Red}$-$N$-injective modules is again
        $\text{Red}$-$N$-injective. In particular, a finite direct sum of
        $\text{Red}$-injective (resp., strongly $\text{Red}$-injective)
        modules is again $\text{Red}$-injective           (resp., strongly $\text{Red}$-injective);
  \item a direct summand of $\text{Red}$-quasi-injective (resp., $\text{Red}$-injective,
        strongly $\text{Red}$-injective) module is again $\text{Red}$-quasi-injective
        (resp., $\text{Red}$-injective, strongly $\text{Red}$-injective) module.
\end{enumerate} 
\end{corollary}

\begin{proposition}
 If  $M$ is  a Noetherian $R$-module, then  a direct sum of $\text{Red}$-$M$-injective
 modules is $\text{Red}$-$M$-injective. 
\end{proposition}

\begin{prf}
  For $D=\bigoplus_{i\in I}D_{i}$, a direct sum of
  $\text{Red}$-$M$-injective modules, let $f:K\rightarrow D$ be
  an $R$-homomorphism, where $K$ is any semi-reduced submodule of $M$.
  Since $K$ is finitely generated, $f(K)\leq \bigoplus_{i=1}^{n}D_{i}$
  for some positive integer $n$.
  Since $\bigoplus_{i=1}^{n}D_{i}$ is $\text{Red}$-$M$-injective,
  then $f$ can be extended to an $R$-homomorphism $\hat{f}:M\rightarrow D$.
\end{prf}

\begin{corollary}
 Let $R_{R}$ be a  Noetherian module. Then, a direct sum of $\text{Red}$-injective
 modules is $\text{Red}$-injective. 
\end{corollary}

\begin{proposition}
        Let $\{M_{i}:i\in I\}$ be a family of $R$-modules and $N$ be an
        $R$-module. Then, $N$ is $\text{Red}$-$(\bigoplus_{i\in I}M_{i})$-injective
        if and only if it is $\text{Red}$-$M_{i}$-injective for each $i$. 
\end{proposition}

\begin{prf}
\begin{description}
   \item[$(\Rightarrow).$] Suppose that $N$ is
        $\text{Red}$-$(\bigoplus_{i\in I}M_{i})$-injective.
        Let $f:\text{Red}(M_{i})\rightarrow N$ be any $R$-homomorphism. By
        hypothesis, any $R$-homomorphism
        $g:\text{Red}(\bigoplus_{i\in I}M_{i})\rightarrow N$ extends to
        $\bar{g}:\bigoplus_{i\in I}M_{i}\rightarrow N$. The required extension
        of $f$ is $\bar{g}\circ \iota$ where $\iota$ is the injection
        $\iota: M_{i}\rightarrow \bigoplus_{i\in I}M_{i}$.
    \item[$(\Leftarrow).$] Suppose that $N$ is $\text{Red}$-$M_{i}$-injective for
        each $i\in I$.
        Since $N$ is $\text{Red}$-$M_{i}$-injective for each $i\in I$;
        let $\theta_{i}:M_{i}\rightarrow N$ be the extension of
        $f_{i}:\text{Red}(M_{i})\rightarrow N$ for each $i\in I$.
        Let also $g:\text{Red}(\bigoplus_{i\in I}M_{i})\rightarrow N$ be any
        $R$-homomorphism. By the fundamental property of direct sum of modules,
        there exists an $R$-homomorphism
        $\theta = \langle \theta_{i}\rangle :\bigoplus_{i\in I}M_{i}\rightarrow N$
        such that $\theta \circ \iota_{i}=\theta_{i}$ for all $i\in I$;
        where $\iota_{i}:M_{i}\rightarrow \bigoplus_{i\in I}M_{i}$ is the
        injection $R$-homomorphism for each $i\in I$. Then $\theta$ is an extension
        of $g:\text{Red}(\bigoplus_{i\in I}M_{i})\rightarrow N$.
        Hence, $N$ is $\text{Red}$-$(\bigoplus_{i\in I}M_{i})$-injective.
\end{description}
 \end{prf}

 \begin{corollary}
 If $A$, $B$, $C$, and $Q$ are $R$-modules and the short exact
 sequence $\{0\} \rightarrow A\xrightarrow{\mu} B\xrightarrow{\varepsilon} C\rightarrow \{0\}$
 splits, then the following conditions hold:
\begin{enumerate}
  \item $Q$ is $\text{Red}$-$C$-injective if and only if it is $\text{Red}$-$(B/\mu (A))$-injective.
  \item $Q$ is $\text{Red}$-$B$-injective if and only if it is $\text{Red}$-$A$-injective
        and $\text{Red}$-$C$-injective.
\end{enumerate} 
\end{corollary}
\begin{prf}
\begin{enumerate}
  \item This follows from the fact that $B/\mu (A)\cong C$.
  \item Follows from Proposition 4, Theorem 1 and the fact that $B\cong A\oplus C$.
\end{enumerate}
\end{prf}

 \begin{proposition}
  Let $N$ and $M$ be $R$-modules. Then the following conditions hold:
\begin{enumerate}
  \item $M$ is injective $\Rightarrow$ $M$ is $N$-injective $\Rightarrow$ $M$
      is $\text{Red}$-$N$-injective $\Rightarrow$ $M$ is $\text{Soc}$-$N$-injective
      $\Rightarrow$ $M$ is min-$N$-injective.
  \item $M$ is injective $\Rightarrow$ $M$ is strongly $\text{Red}$-injective
      $\Rightarrow$ $M$ is strongly $\text{Soc}$-injective $\Rightarrow$ $M$ is strongly
      $min$-injective $\Leftrightarrow$ $M$ is strongly simple-injective.
\end{enumerate} 
\end{proposition}

\begin{prf}
  Elementary.
\end{prf}

 \begin{proposition}
  For an $R$-module $M$, if $\text{Red}(M)$ is a direct summand of $M$,
 then every $R$-module is $\text{Red}$-$M$-injective. 
\end{proposition}
\begin{prf}
  Suppose that $K$ is an $R$-module and $\text{Red}(M)\subseteq^{\oplus} M$.
  We show that $K$ is $\text{Red}$-$M$-injective.
  Let $f:\text{Red}(M)\rightarrow K$ be any $R$-homomorphism.
  Since $\text{Red}(M)$ is a direct summand of $M$, there exists a proper $R$-submodule $P$
  of $M$ such that $M=\text{Red}(M)\oplus P$.
  There exists an $R$-homomorphism $f':M\rightarrow \text{Red}(M)$
  such that $f'(n+p)=n$, for all $n\in \text{Red}(M)$ and $p\in P$.
  Then, the $R$-homomorphism $f\circ f':M\rightarrow K$
  is an extension of $f$ because $(f\circ f')(n+p)=f(f'(n+p))=f(n)$ for all $n+p~\in M$.
  Hence, $K$ is $\text{Red}$-$M$-injective.
 \end{prf}

 \begin{theorem}
        For a projective $R$-module $M$, the following conditions are equivalent:
\begin{enumerate}
  \item Every quotient of a $\text{Red}$-$M$-injective $R$-module is $\text{Red}$-$M$-injective.
  \item Every quotient of an injective $R$-module is $\text{Red}$-$M$-injective.
  \item $\text{Red}(M)$ is a projective $R$-module.
\end{enumerate} 
\end{theorem}

\begin{prf}
\begin{description}
 \item[$(1\Rightarrow 2).$] This is due to the fact that every injective $R$-module
 is $\text{Red}$-$M$-injective.

 \item[$(2\Rightarrow 3).$] Consider the diagram in Figure 3 below:
\begin{figure}[H]
  \centering
  \input kitare.pic
  \label{kitare}
  \caption{}
\end{figure}

  where $E$ and $N$ are $R$-modules, $\varepsilon$ an $R$-epimorphism,
  and $f$ an $R$-homomorphism. By \cite[Proposition 5.1]{eilenberg1956homological},
  assume that $E$ is injective. Since $N$ is $\text{Red}$-$M$-injective $f$ can
  be extended to an $R$-homomorphism $g:M\rightarrow N$. Since $M$ is projective,
  $g$ can be lifted to an $R$-homomorphism $\widetilde{g}:M\rightarrow E$
  such that $\varepsilon \circ \widetilde{g} =g$.
  Define $\widetilde{f}:\text{Red}(M)\rightarrow E$ by
  $\widetilde{f}=\widetilde{g}|_{_{\text{Red}(M)}}$.
  Then $\varepsilon \circ \widetilde{f}=
  \varepsilon \circ \widetilde{g}|_{_{\text{Red}(M)}}=f$.
  Hence, $\text{Red}(M)$ is projective.
  \item[$(3\Rightarrow 1).$] Let $N$ and $L$ be $R$-modules with
  $\varepsilon :N\rightarrow L$
  an $R$-epimorphism and $N$ is $\text{Red}$-$M$-injective.
  Consider the diagram in Figure 4.
\begin{figure}[H]
  \centering
  \input keza.pic
  \label{keza}
  \caption{}
\end{figure}

  Since $\text{Red}(M)$ is projective, $f$ can be lifted to an $R$-homomorphism
  $g:\text{Red}(M)\rightarrow N$ such that
  $\varepsilon \circ g(m)=f(m)$, for all $m\in \text{Red}(M)$.
  Since $N$ is $\text{Red}$-$M$-injective, $g$ can be extended to
  an $R$-homomorphism $\widetilde{g}:M\rightarrow N$.
  Hence, $\varepsilon \circ \widetilde{g}:M\rightarrow L$ extends $f$.
  \end{description}
\end{prf}

\begin{corollary}
  The following conditions are equivalent for a reduced projective $R$-module:
\begin{enumerate}
  \item Every quotient of a $\text{Red}$-injective $R$-module is $\text{Red}$-injective.
  \item Every quotient of an injective $R$-module is $\text{Red}$-injective.
  \item $\text{Red}(R_{R})$ is a projective module.
\end{enumerate}
 In addition, if every semi-reduced submodule of a projective $R$-module is projective,
 then $\text{Red}(R_{R})$ is a projective module. 
 \end{corollary}
 
\begin{prf}
  $1\Leftrightarrow 1\Leftrightarrow4$ follows from Theorem 2.
  The additional case follows from the fact that $\text{Red}(R_{R})$ is a semi-reduced
  submodule of a projective module $R_{R}$.
\end{prf}

\begin{proposition}
   Let $R$ be a Principal Ideal Domain (PID) and $N$ be an $R$-module.
  Then, the following statements hold:
\begin{enumerate}
  \item If every free $R$-module is $\text{Red}$-$N$-injective
        then each of its submodules is $\text{Red}$-$N$-injective.
  \item If every projective $R$-module is $\text{Red}$-$N$-injective
        then each of its submodules is $\text{Red}$-$N$-injective.
  \item Every projective $R$-module is $\text{Red}$-$N$-injective if
        and only if every free $R$-module is $\text{Red}$-$N$-injective.
\end{enumerate} 
\end{proposition}

\begin{prf}
\begin{enumerate}
  \item Suppose that every free $R$-module $M$ is $\text{Red}$-$N$-injective,
        and $L\leq M$. Since over
        a PID a submodule of a free module is free, $L$ is free.
        By hypothesis, $L$ is $\text{Red}$-$N$-injective.
  \item Suppose that every projective $R$-module $P$ is $\text{Red}$-$N$-injective,
        and $K\leq P$. Since over a PID
        a submodule of a projective $R$-module is projective, $K$ is projective.
        By hypothesis, $K$ is $\text{Red}$-$N$-injective.
  \item  Over a PID every projective module is free.
         The converse holds since any free module is projective.
\end{enumerate}
\end{prf}

\begin{definition}
 Let $X$ be a submodule of a module $M$. We say that $\text{Red}(M)$ respects
 $X$ if there exists a direct summand $A$ of $M$ contained in $X$
 such that $X = A\oplus B$ and $B \leq \text{Red}(M)$. $M$ is called
 $\text{Red}(M)$-lifting if $\text{Red}(M)$ respects every submodule of $M$.
\end{definition}

\begin{proposition}
 Let $N$ be an $R$-module. If $N$ is $\text{Red}(N)$-lifting,
 then any $R$-module $K$ is $\text{Red}$-$N$-injective
 if and only if $K$ is $N$-injective. 
\end{proposition}

\begin{prf}
\begin{description}
  \item[$(\Rightarrow).$] Suppose that $K$ is $\text{Red}$-$N$-injective.
     Let $L$ be any submodule of $N$, $\iota :L\rightarrow N$
     the inclusion map and $f:L\rightarrow K$ any $R$-homomorphism.
     Since $\text{Red}(N)$ respects $L$, $L$ has a decomposition $L=A\oplus B$
     such that $A\subseteq^{\oplus} N$ and $B\leq \text{Red}(N)$. $N=A\oplus A'$
     for some submodule $A'$ of $N$. Then, $L=A\oplus (L\cap A')$
     and $L\cap A'$ is semi-reduced. Let $i: L\cap A'\rightarrow L$
     be the inclusion map and $f|_{L\cap A'}:L\cap A'\rightarrow K$.
     Since $K$ is $\text{Red}$-$N$-injective, there exists an $R$-homomorphism
     $g:N\rightarrow K$ such that $g\circ \iota \circ i=f|_{L\cap A'}$.
     Now, define $h:N\rightarrow K$ by $h(a+a')=f(a)+g(a')$~~~$(a\in A, a'\in A')$.
     Then $h\circ \iota =f$, and hence $K$ is $N$-injective.
  \item[$(\Leftarrow).$] Every $N$-injective module is $\text{Red}$-$N$-injective.
      This is due to the fact that for every $N$-injective module $K$,
      any $R$-homomorphism from any submodule of $N$ to $K$ extends to $N$.

\end{description}
\end{prf}

  Note that a semi-simple module as well as a module $M$ with no reduced submodule
  (i.e., one for which $\text{Red}(M)=\{0\}$) is $\text{Red}(M)$-lifting.

 Let $K$ and $L$ be submodules of $M$. $M$ is said to satisfy:
 \begin{enumerate}
 \item C1-condition if every submodule of $M$ is essential in a summand of $M$.
 \item C2-condition if $K \cong L$ and $K\subseteq^{\oplus} M$,
          then $L\subseteq^{\oplus} M$.
 \item C3-condition if $K\cap L =\{0\}$, $K\subseteq^{\oplus} M$
          and $L\subseteq^{\oplus} M$, then $K\oplus L\subseteq^{\oplus} M$.
 \end{enumerate}

 A module is {\it quasi-continuous} if it satisfies C1 and C3 conditions.

 Proposition 9 shows that $\text{Red}$-quasi-injective modules inherit
 a weaker version of C2-condition and C3-conditions.

 \begin{proposition}
  Suppose that an $R$-module $N$ is $\text{Red}$-quasi-injective.
 \begin{enumerate}
   \item ($\text{Red}$-C2) If $P$ and $Q$ are semi-reduced submodules of $N$,
                            $P\cong Q$ and $P\subseteq^{\oplus} N$,
                            then $Q\subseteq^{\oplus} N$.
   \item ($\text{Red}$-C3) Let $P$ and $Q$ be semi-reduced submodules of $N$
                            with $P\cap Q=\{0\}$. If $P\subseteq^{\oplus}N$
                            and $Q\subseteq^{\oplus}N$;
                            then $P\bigoplus Q \subseteq^{\oplus}N$.
 \end{enumerate} 
 \end{proposition}

 \begin{prf}
 \begin{enumerate}
   \item Since $P\cong Q$, and $P$ is $\text{Red}$-$N$-injective, being a direct
     summand of a $\text{Red}$-quasi-injective module $N$, $Q$ is
     $\text{Red}$-$N$-injective by Corollary 2(2).
     If $i:Q\rightarrow N$ is the inclusion map, the identity
     $id_{Q}:Q\rightarrow Q$ has an extension $\eta:N\rightarrow Q$
     such that $\eta \circ i=id_{Q}$, and hence $Q\subseteq^{\oplus} N$.
   \item Since both $P$ and $Q$ are direct summands of $N$; then both $P$
         and $Q$ are $\text{Red}$-$N$-injective. Then the semi-reduced
         module $P\oplus Q$ is $\text{Red}$-$N$-injective, and so a
         direct summand of $N$ by an argument similar to the one given in 1.
 \end{enumerate}
\end{prf}

\section{Strongly \text{Red}-injective modules}
 In this section, we characterize quasi-Frobenius rings and right $V$-rings
 in terms of strongly $\text{Red}$- injective modules.
 A ring $R$ is called {\it right semi-Artinian} if every
 non-zero $R$-module has nonzero socle. A submodule $S\leq M$ is small if,
 for any submodule $N\leq M$, $S+N=M$ implies that $N=M$. The projective cover
 of an $R$-module $M$ is a projective module $P$ for which there is an
 epimorphism $P\rightarrow M$ whose kernel is small. A ring
 $R$ is left perfect if every $R$-module has a projective cover.

 \begin{proposition}
 The following implications hold: 

   $R$ is right semi-Artinian $\Rightarrow$
  every strongly Red-injective $R$-module is injective $\Rightarrow$
   every strongly Red-injective $R$-module is quasi-continuous. 

 In particular, over a left perfect ring $R$, every strongly $\text{Red}$-injective
 right $R$-module is injective. 
\end{proposition}

\begin{prf}
   For a right semi-Artinian ring $R$, suppose that a non-zero $R$-module $M$ is strongly
   $\text{Red}$-injective.
   Then, $\{0\}\neq \text{Soc}(M)\subseteq^{e}M$.
   Amin  \emph{et al}., in \cite[Corollary 3.2]{amin2005} showed that a strongly
   $\text{Soc}$-injective module with essential socle is injective.
   Since $M$ has essential socle, it is injective.
   $M$ is quasi-continuous because every injective module is
   quasi-continuous see \cite[p.18]{mohamed1990continuous}.
   The last statement follows from the fact that every left perfect ring is right
   semi-Artinian, see \cite[Theorem 11.6.3]{kasch1982modules}.
\end{prf}

 A ring $R$ is called {\it quasi-Frobenius} if $R$ is right (or left) Artinian and
 right (or left) self-injective. Equivalently, $R$ is {\it quasi-Frobenius}
 if and only if every injective $R$-module is projective if and only if
 every projective $R$-module is injective.
 A ring whose all simple right modules are injective is called a right {\it $V$-ring}.

\begin{theorem}
 A ring $R$ is quasi-Frobenius if and only if every strongly $\text{Red}$-injective
 module is projective. 
\end{theorem}

\begin{prf}
 If $R$ is quasi-Frobenius, then $R$ is right semi-Artinian and so by Proposition 10 every strongly
 $\text{Red}$-injective module is injective, and hence projective since $R$ is quasi-Frobenius.
 Conversely, if every strongly $\text{Red}$-injective module is projective,
 then in particular every injective module is projective, and so $R$ is quasi-Frobenius.
\end{prf}

 \begin{theorem}
  $R$ is a right $V$-ring if and only if every simple $R$-module is strongly
 $\text{Red}$-injective. 
\end{theorem}

\begin{prf}
 Suppose that $M$ is a simple $R$-module where $R$ is a right $V$-ring.
 Then, by definition of a $V$-ring, $M$ is injective.
 Hence, $M$ is strongly $\text{Red}$-injective.  Conversely, suppose that
 any simple module $M$ is strongly $\text{Red}$-injective.
 Since $M$ is simple, $\text{Soc}(M)=M$ and hence
 $\{0\}\neq \text{Soc}(M)\subseteq^{e}M$. Since $M$ has essential socle,
 it is injective by \cite[Corollary 3.2]{amin2005}. Hence, $R$ is a right $V$-ring.
\end{prf}

\begin{corollary}
 Let $M$ be an $R$-module with essential socle. The following statements are
 equivalent:
 \begin{enumerate}
   \item $M$ is injective.
   \item $M$ is strongly $\text{Red}$-injective.
   \item $M$ is strongly $\text{Soc}$-injective.
 \end{enumerate} 
 
\end{corollary}

\begin{prf}
  By Proposition 5(2), $1\Rightarrow 2\Rightarrow 3$. By 
   \cite[Corollary 3.2]{amin2005}, $3\Rightarrow 1$ which completes the proof.
\end{prf}

\begin{remark}
   We make the following observations:
\begin{enumerate}
 \item It is easy to check that any sort of injectivity that lies between injective and
  strongly $\text{Red}$-injective modules would lead to Theorems 3 and 4.
  \item $\text{Red}$-injectivity is a less restricted notion than injectivity but carries
  most of the properties of injectivity.
 \item $\text{Red}$-injectivity is  much closer to injectivity than $\text{Soc}$-injectivity.
 
 \item When the ring is not commutative,  a semi-simple module  need not be semi-reduced, see Example 3.1 below:
  \end{enumerate}
\end{remark}

 \begin{example}\label{ex}
  Let the ring $R$ be the collection of all $2\times 2$ matrices over
  the field of real numbers. The module $M=R_{R}$ is semi-simple but
  not reduced. For if $m=\left(
                           \begin{array}{cc}
                             1 & 0 \\
                             0 & 1 \\
                           \end{array}
                         \right)
  \in M$ and $r=\left(
                  \begin{array}{cc}
                    1 & -1 \\
                    1 & -1 \\
                  \end{array}
                \right)
  \in R$, then $mr\neq 0$ but $mr^{2}=0$.
 Since a direct sum of reduced modules is reduced,
 $M$ is a direct sum of simple modules which is not reduced. A simple module
 over a not necessarily commutative ring need not be reduced. $M$ is not
 semi-reduced.
 \end{example}

  The following implications hold:
  
  Injective $\Rightarrow$ quasi-injective $\Rightarrow$ pseudo-injective $\Rightarrow$
  $\text{Red}$-quasi-injective.
  
  For the first two implications, see \cite{singh}. The last  implication is trivial, it follows directly from the definitions. 
   
  \begin{example}
  Let $R$ be the ring of all eventually constant sequences $(x_{n})_{n\in \mathbb{N}}$
  of elements in $\mathbb{F}_{2}$, the field of two elements.
  Then, $E(R_{R})=\prod_{n\in \mathbb{N}}\mathbb{F}_{2}$, which has only one automorphism,
  namely the identity automorphism. By \cite[Example 9]{er2013rings}, $R_{R}$
  is pseudo-injective but it is not quasi-injective. It therefore follows that $R_{R}$ is
  $\text{Red}$-quasi-injective but not injective.
  \end{example}

  An $R$-module $M$ is said to satisfy the {\it exchange property} if for every $R$-module
  $A$ and any two direct sum decomposition $A=M'\bigoplus N=\bigoplus_{i\in I}A_{i}$ with
  $M'=M$, there exists a submodule $B_{i}$ of $A_{i}$ such that
  $A=M'\bigoplus (\bigoplus_{i\in I}B_{i})$. An $R$-module is called {\it clean} if
  its endomorphism ring, $\text{End}_{R}(M)$ is clean, i.e, for all $f\in \text{End}_{R}(M)$,
  $f=e+u$ with $e$ idempotent and $u$ a unit. Pseudo-injective modules
  (and hence quasi-injective and injective modules) are clean and also satisfy the
  exchange property, see \cite{asensio2015automorphism} and \cite{asensio2013automorphism}. Note that pseudo-injective modules are equivalent to
   automorphism-invariant modules as they are being referred to in  \cite{asensio2015automorphism} and \cite{asensio2013automorphism}. The equivalence
    was proved in  \cite{er2013rings}.     We now ask:
    
  \begin{question}
  Are the $\text{Red}$-quasi-injective modules and $\text{Soc}$-quasi-injective modules clean?
  Do they satisfy the exchange property?
  \end{question}

\end{document}